\documentclass[oneside,reqno, twoside, 11pt]{amsart}
\usepackage{color}
\usepackage[T1]{fontenc}
\usepackage{lmodern}
\usepackage[english]{babel}
\usepackage{mathrsfs}  
\usepackage{upgreek}
\usepackage{amsfonts}
\usepackage{extpfeil}
\usepackage{verbatim}
\usepackage{BOONDOX-cal}
\usepackage{microtype}
\usepackage[new]{old-arrows}
\usepackage{amsmath}
\usepackage[all, cmtip]{xy}
\usepackage{tikz-cd}

\usepackage{amsthm}
\usepackage{mathtools}

\usepackage{bbm}
\usepackage{paralist}
\usepackage[T1]{fontenc}
\usepackage[colorlinks,citecolor=magenta,pagebackref=true,urlcolor=magenta, bookmarksdepth=3]{hyperref}
\usepackage[nodisplayskipstretch]{setspace}
\usepackage{breqn}

\usepackage{nicefrac}
\usepackage[font=small,labelfont=bf]{caption}
\usepackage[labelformat=simple]{subcaption}

\usepackage{epsfig}

\DeclareFontFamily{OT1}{pzc}{}
\DeclareFontShape{OT1}{pzc}{m}{it}{<-> s * [1.10] pzcmi7t}{}
\DeclareMathAlphabet{\mathpzc}{OT1}{pzc}{m}{it}

\makeatletter
\newtheorem*{rep@theorem}{\rep@title}
\newcommand{\newreptheorem}[2]{%
	\newenvironment{rep#1}[1]{%
		\def\rep@title{#2~\ref{##1}}%
		\begin{rep@theorem}}%
		{\end{rep@theorem}}}

\makeatother
\theoremstyle{plain}
\newtheorem*{thm*}{Theorem}
\newtheorem{thm}{Theorem}[section]
\newreptheorem{mthm}{Main Theorem}
\newreptheorem{mcor}{Corollary}
\newreptheorem{mlem}{Lemma}
\newreptheorem{thm}{Theorem}

\newtheorem{cor}[thm]{Corollary}
\newtheorem{lem}[thm]{Lemma}

\newtheorem*{lem*}{Lemma}

\newtheorem{conj}[thm]{Conjecture}

\newtheorem{qst}[thm]{Question}
\newtheorem{prb}[thm]{Problem}

\theoremstyle{definition}

\newcommand{\supp}{\mathrm{supp}\hspace{1mm}}

 \newcommand{\R}{\mathbb{R}}

\newcommand{\RR}{\Psi}

\newcommand{\tet}{\theta}

\newcommand{\st}{\mr{st}}

\newcommand{\sbseq}{\subseteq}
\newcommand{\spseq}{\supseteq}

\newcommand{\vanish}[1]{}

\def\V{{\bf V}}

\def\im{\mathrm{im}}%\hspace{1mm}}
\def\ker{\mathrm{ker}}%\hspace{1mm}}

\def\sbs\subset
\def\sbseq{\subseteq}

\def\langle{\left<}
\def\rangle{\right>}

\def\({\left(}
\def\){\right)}
\def\no={\,{\,|\!\!\!\!\!=\,\,}}

\def\no={\,{\,|\!\!\!\!\!=\,\,}}
\def\sbseq{\subseteq}

\def\sbseq{\subseteq}
\def\sbs\subset
\def\spseq{\supseteq}

\newcommand{\xqedhere}[2]{%
	\rlap{\hbox to#1{\hfil\llap{\ensuremath{#2}}}}}

\newcommand\Defn[1]{\textbf{#1}}
\newcommand{\cm}[1]{}
\newcommand\mc[1]{\mathcal{#1}}

\newcommand\mbf[1]{\mathbf{#1}}

\newcommand\mr[1]{\mathrm{#1}}

\newcommand{\bigslant}[2]{{\raisebox{.3em}{$#1$} \Big/ \raisebox{-.3em}{$#2$}}}

\newcommand\x{\mathbf{x}}

\DeclareMathOperator{\lk}{lk}

\DeclareMathOperator{\cone}{cone}

\DeclareMathOperator{\Lk}{lk}

\DeclareMathOperator{\St}{st}

\title{FAQ on the g-theorem and the hard Lefschetz theorem for face rings}
\author{Karim Adiprasito }

\address{\emph{Karim Adiprasito}, Einstein Institute of Mathematics, Hebrew University of Jerusalem, Jerusalem, Israel \emph{and} Department of Mathematics, KTH Stockholm, Stockholm, Sweden}
\email{adiprasito@math.huji.ac.il}

\date{24.06.2019}

%Primary  14M25, 05C38 ; Secondary 32S50, 52C25,  13F55}

\begin{document}
	
	\begin{abstract}
We review the hard Lefschetz theorem for simplicial spheres, as well as the theory at its core: perturbations of maps, biased Poincar\'e pairings and a cobordism argument that relates the Lefschetz property of a manifold to the Lefschetz property on its boundary. We also sketch an alternative argument based on edge-contractions.
	\end{abstract}
	
	\maketitle
	
	\newcommand{\AR}{\mathcal{A}}
	\newcommand{\BR}{\mathcal{B}}
	\newcommand{\CR}{\mathcal{C}}
	\newcommand{\Mu}{M}
	\newcommand{\Soc}{\mathcal{S}\hspace{-1mm}\mathcal{o}\hspace{-1mm}\mathcal{c}}
	\newcommand{\Socl}{{\Soc^\circ}}

\section{Counting}
	
One of the most basic problems that almost every combinatorialist asks is to count things. And when you ask a combinatorial topologist, they usually want to count things associated to topological spaces. 

More precisely: many topological spaces can be triangulated, that is, there are simplicial complexes that encode them. A tetrahedron, for instance, triangulates a $3$-dimensional ball. The boundary of that tetrahedron triangulates the $2$-dimensional sphere. But so does the boundary of the icosahedron. Now we can count the number of faces that each of these triangulations have:

In the former case, we get one empty face, four vertices, six edges, and four triangles. Nice and symmetric (well, almost). This is usually recorded in a vector called the $f$-vector recording the number of $i$-dimensional faces $f_i$, which reads in this case as
\[(1,4,6,4).\]

In the latter, we get one empty face, twelve vertices, thirty edges, and twenty triangles (where it gets its name), and hence an $f$-vector
\[(1,12,30,20).\]
This is no longer so nice and symmetric. However, Sommerville \cite{Sommerville} had the idea to define another vector
\[h_k \ :=\  \sum_{i=0}^k (-1)^{k-i}\binom{d-i}{k-i}f_{i-1}.\]
And now, there is a small miracle: he then established\footnote{He stated this for boundaries of simplicial polytopes, but the proof works for simplicial spheres} that for a simplicial sphere of dimension $(d-1)$, we have
\[h_k \ =\  h_{d-k}.\]
These are the so-called Dehn-Sommerville relations. This means in particular that there are  nontrivial linear relations between the face numbers of simplicial polytopes, and that everything is defined from the first half of the entries. So, McMullen \cite{zbMATH03333960} had the idea to consider another vector:
\[g_k\ := h_k-h_{k-1}\ \text{for}\ k\le \frac{d}{2}\]
He formulated the following conjecture
\begin{conj}
	A vector of $d$ integers is the $f$-vector of a simplicial sphere $\Sigma$ if and only if the associated $g$-vector is an $M$-sequence, that is, there is a quotient $Q$ of a polynomial ring $\R[\x]$ by a homogeneous ideal so that \[g_i(\Sigma)\ =\ \dim Q^i.\]
\end{conj}
At this point, Stanley entered. He realized that there is at least always a ring that encodes the $h$-vector.
\section{Rings}
If $\Delta$ is an abstract simplicial complex defined on the groundset $[n]\coloneqq \{1,\cdots,n\}$, let $I_\Delta\coloneqq \langle \x^{\mbf{a}}:\ \supp(\mbf{a})\notin\Delta\rangle$ denote the nonface ideal in $\R[\x]$, where $\R[\x]=\R[x_1,\cdots,x_n]$. Let $\R^\ast[\Delta]\coloneqq \R[\x]/I_\Delta$ denote the \Defn{face ring} of $\Delta$. Now, we pick a sufficient number of linear forms to make sure the quotient is finite dimensional:

We may associate to the vertices of $\Delta$ the coordinates $\V_\Delta=({v_1},\cdots, {v_n}) \in \R^{l\times n}$, obtaining a system of linear forms by considering  $\V_\Delta\x = \Theta$. With this, we obtain a \Defn{geometric simplicial complex}.

The face ring of a {geometric simplicial complex} $\Delta$ is considered with respect to its natural system of parameters induced by the coordinates, that is, 
\[\AR^\ast(\Delta)\ :=\ \R^\ast[\Delta]/\Theta\R^\ast[\Delta].\]

A geometric simplicial complex in $\R^d$ is \Defn{proper} if the image of every $k$-face, with $k<d$, linearly spans a subspace of dimension $k+1$. If $\Delta$ is of dimension $(d-1)$, and is given a proper coordinates in $\R^d$, then $\AR^\ast(\Delta)$ is finite-dimensional as a vector space. 

Stanley observed \cite{StanleyCMC}, based on a theorem of Reisner \cite{Reisner}:

\begin{thm}
For a triangulated sphere $\Sigma$ of dimension $(d-1)$, 
\[h_i(\Sigma)\ =\ \dim \AR^i(\Sigma).\]
\end{thm}

\section{Stanley and Lefschetz} Here, Stanley observed \cite{StanleyHL} that McMullen's conjecture is true if there existed an $\ell$ in $\AR^1[\Sigma]$ so that
\[\AR^{i}(\Sigma)\ \xrightarrow{\ \cdot\ell\ }\ \AR^{i+1}(\Sigma)\]
is injective for $i\le \frac{d}{2}-1$, or stronger if 
\[\AR^k(\Sigma)\ \xrightarrow{\ \cdot \ell^{d-2k} \ }\ \AR^{d-k}(\Sigma). \]
is an isomorphism for every $k\le \frac{d}{2}$. 
The former is known as the \Defn{weak Lefschetz property}, the latter as the \Defn{hard Lefschetz property}. He needed this for \emph{some} geometric realization of the simplicial complex. 

And amazingly, Stanley then observed that the hard Lefschetz property is actually true for spheres that arise as boundaries of simplicial polytopes (with respect to their given geometric realization) using deep results in algebraic geometry.
Specifically, if $\Sigma$ is realized as the boundary of a polytope, then the class of a convex function acts as the desired Lefschetz element (we refer to \cite{KF} for an excellent exposition). But the general case remained open, and only recently, the hard Lefschetz theorem was proven in this generality. It is useful to state this in a more general perspective for the purposes of this paper: 

 A \Defn{relative simplicial complex} $\RR=(\Delta, \Gamma)$ is a pair of simplicial complexes $\Delta, \Gamma$ with $\Gamma \subset \Delta$.
If $\RR=(\Delta,\Gamma)$ is a relative simplicial complex, then we can define the \Defn{relative face module}
\[\R^\ast[\RR]\ \coloneqq \ \bigslant{I_\Gamma}{I_\Delta}.
\]
The relative face ring of a geometric complex is then defined as
\[\AR^\ast(\RR)\ :=\ \R^\ast[\RR]/\Theta\R^\ast[\RR].\]

\begin{thm}\label{mthm:gl}
		Consider a triangulated $(d-1)$-sphere or ball $\Delta$, and the associated graded commutative ring $\R[\Delta]$. Then
		there exists an open dense subset of the Artinian reductions $\mathcal{R}$ of $\R[\Delta]$ and an open dense subset $\mathcal{L} \subset \AR^1(\Delta)$, where $\AR(\Delta)\in \mc{R}$, such that for every $k\le\nicefrac{d}{2}$, we have the
{hard Lefschetz property} for every $\AR(\Delta)\in \mathcal{R}$ and every $\ell \in \mathcal{L}$:
			\[\AR^k(\Delta,\partial \Delta)\ \xrightarrow{\ \cdot \ell^{d-2k} \ }\ \AR^{d-k}(\Delta)\]
			is an isomorphism.
\end{thm}

There are two proofs of this theorem: one is detailed in \cite{AHL}, the other, slightly different and relying on edge contractions, is sketched here. However, their core, the theory of biased pairings, is the same, and is what will be presented here. We want to give a brief overview over the key ideas of these proofs. Up until the final two sections, they are essentially the same, with some variations in generality.

Let us first recall a basic ingredient to the inductive structure of the proofs: Two classical lemmas that help tremendously in reducing from higher dimensional cases to lower-dimensional ones.

\section{Two cone lemmas} Recall that the  \Defn{star} and \Defn{link} of a face $\sigma$ in $\Delta$ are
the subcomplexes \[\St_\sigma \Delta\ \coloneqq \ \{\tau:\exists \tau'\supset \tau,\ \sigma\subset
\tau'\in \Delta\}\ \
\text{and}\ \ \Lk_\sigma \Delta\ \coloneqq \ \{\tau\setminus \sigma: \sigma\subset
\tau\in \Delta\}.\]
For geometric simplicial complexes $\Delta$, we shall think of the star of a face as a geometric subcomplex of $\Delta$, and the link of a face $\sigma$ as the geometric simplicial complex obtained by the orthogonal projection to $\mr{span}(\sigma)^\bot.$
Let us denote the \Defn{deletion} of $\sigma$ by $\Delta-\sigma$, the maximal subcomplex of $\Delta$ that does not contain $\sigma$.
Let 
\[\St_\sigma^\circ \Delta\ \coloneqq\  (\St_\sigma \Delta,\St_\sigma \Delta-\sigma).\]
We have the following two elementary lemmas.

\begin{lem}[Cone lemma I, see {\cite[Thm. 7]{Lee}}]
	For any vertex $v\in \Delta$, where $\Delta$ is a geometric simplicial complex in $\R^d$, and for any integer $k$, we have an isomorphism
	\[\AR^k({\Lk}_v \Delta)\ \cong\ \AR^{k}(\St_v \Delta).\]
\end{lem}

\begin{lem}[Cone lemma II, see {\cite[Lem. 3.3]{AdiprasitoTC}}]
	In the situation of the first cone lemma we have a natural isomorphism
	\[\AR^{k}(\St_v \Delta)\ \xrightarrow{\ \cdot x_v \ }\ \AR^{k+1}(\St_v^\circ \Delta).\]
\end{lem}

We come to the first core principle, a key new idea for the proof of the hard Lefschetz property.

\section{Generic combinations of linear maps}

First is the idea to construct the map $\ell$ iteratively. We rely on the following principle:

\begin{lem}[Lemma 6.1, \cite{AHL}]\label{lem:perturbation}
	Given two linear maps 
	\[A, B: \mathcal{X}\ \longrightarrow\ \mathcal{Y}\]
	of two vector spaces $\mathcal{X}$ and $\mathcal{Y}$ over $\R$.
	Assume that 
		\[B(\ker A)\ \cap\ \im A\ =\ {0}\ \subset\ \mathcal{Y} .\]
		Then a generic linear combination 
		$A ``{+}" B$ of $A$ and $B$
		has kernel 
		\[\ker (
		A\ ``{+}"\ B)\ = \ \ker A\ \cap\ \ker B.\]
\end{lem}

This lemma may feel unnatural, but it is exactly fitting for the spaces we consider. As it turns out, Gr\"abe \cite{Grabe} showed that there is a perfect pairing, the Poincar\'e pairing,
\[\AR^k(\Delta,\partial \Delta)\ \times\ \AR^{d-k}(\Delta)\ \longrightarrow\ \AR^{d}(\Delta,\partial \Delta)\ \cong\ \R.\]
Now, we notice that $\ker \cdot \ell^{d-2k}$ and $\im \cdot \ell^{d-2k}$ are orthogonal complements for whatever $\ell$ we choose. Lets assume $d-2k=1$ for now, and look at the map
\[\AR^k(\Delta,\partial \Delta)\ \xrightarrow{\ \cdot \ell\ }\ \AR^{k+1}(\Delta)\]
with $\ell$ to be constructed.

Stay with me: Say now we wish to construct $\ell$ as a Lefschetz map inductively. Then we can start with a vertex $v$ to begin with, and take the associated element $x_v$. This is usually not good enough to induce an isomorphism, but it has controllable kernel and image. In each and every step, we now pick another vertex (say $w$), take the associated element ($x_w$), and generically combine it with what we have (say $\ell'$). Then to make sure that we have an isomorphism at the end, or equivalently an injection, it is good if the kernel of this generic combination is as small as possible.

The ideal case is, if $\ell'$ is supported in vertices $W$, that 
\[\ker \ell'\ =\ \ker ``{\sum_{v\in W}}"\ x_v\ =\ \bigcap_{v\in W} \ker\ x_v\]
This is called the \Defn{transversal prime property} for $W$.

 If we want to prove the analogue for $\ell'\ ``{+}"\ x_w$, then the lemma tells us what to do: show that 
\[B(\ker A)\ \cap\ \im A\ =\ {0}\ \subset\ \mathcal{Y}\]
where $A$ is the multiplication with $\ell'$, and $B$ is the multiplication with $x_w$. 

Now comes the cool part: Multiplication with $x_w$ is the same as pullback to the link of $w$, which is again a disk or a sphere, depending on whether $w$ is in the interior of $\Delta$ or on the boundary. Hence we are trying to show that two orthogonal complements intersect trivially. Which is equivalent to saying that the Poincar\'e pairing is nondegenerate on one (and equivalently both) of these spaces. 

Hence, we want to show that $x_w\ker \ell'$, seen as subspace in $\AR^{\ast}(\lk_w \Delta, \lk_w \partial\Delta)$ via the cone lemmas, has a perfect Poincar\'e pairing induced by the Poincar\'e pairing 
\[\AR^{k}(\lk_w \Delta, \lk_w \partial \Delta) \times \AR^{d-k-1}(\lk_w \Delta)\ \longrightarrow\  \R.\]

This is the notion of biased Poincar\'e duality, the second key part of the proofs. We shall introduce it after an example.

\section{An example}
As it turns out, this property is quite controllable again, by reducing it to Lefschetz isomorphisms. Lets explain how in an example, following \cite[Section 6.5]{AHL}: Assume we are in a sphere $\Sigma$ of dimension $d-1=2k$, and have shown that 

\begin{compactenum}[(1)]
\item There exists a set of vertices $W$ such that $\Sigma-W$ is a disk $\Delta$ and
\item We have the \Defn{transversal prime property} for $W$ in $\Sigma$, that is, the kernel of $\ell' = ``{\sum_{v\in W}}"x_v$ is exactly $\AR^k(\Delta,\partial \Delta),$ the intersection of the kernels of the $x_v, v\in W$.
\end{compactenum}

Lets pick $w$ a vertex in the boundary of $\Delta$. We wish to establish that $\ell'``+" x_w$ has kernel $\AR^k(\Delta-w,\partial (\Delta-w))$. 
Then we wish to understand what to do to ensure that the Poincar\'e pairing of $\AR^k(\st_w\Sigma)$ does not degenerate when restricting to the pullback of $\AR^k(\Delta,\partial \Delta)$. On the other hand,
we have a short exact sequence
\[\AR^k(\st_w \Delta, \st_w \partial \Delta)\ \xrightarrow{\ \cdot x_w\ }\ \AR^{k+1}(\Delta, \Delta-w)\ \longrightarrow\ \AR^{k+1}(\partial \Delta, \partial\Delta-w)\ \longrightarrow\ 0.\]
Hence, the kernel of the multiplication with $x_w$ is
$\AR^k(\Delta-w, \partial \Delta-w)$ provided \[\AR^{k+1}(\partial \Delta, \partial\Delta-w)\ =\ \AR^k({\Lk}_w \partial \Delta)\ =\ 0.\] Where is the Lefschetz theorem though? 

Well
$\uppi$ is a projection of $\R^d/w$ to a hyperplane, and $h$ is the coordinate with respect to that projection, $\tet=h\cdot\x\, $ the associated linear form, then the last condition is equivalent to the Lefschetz property in degree $k-1$: For the map
\[\AR^{k-1}(\uppi{\Lk}_w \partial \Delta)\ \xrightarrow{\ \cdot \tet \ }\ \AR^{k}(\uppi{\Lk}_w \partial \Delta)\]
to be an isomorphism is exactly equivalent to $\AR^k({\Lk}_w \partial \Delta)=0$.
		 
Unfortunately, this does not give a complete proof: We need to make sure that we can remove vertices iteratively so that the deletions remain disks. Also, we need to make sure that the same happens to the ${\Lk}_w \partial \Delta$. This is a rather strong restriction, and introduced as $L$-decomposability in \cite{AHL} (somewhat different from what Billera and Provan introduced originally \cite{BP}, which misses the condition on ${\Lk}_w \partial \Delta$). 

For the sake of clarity, we call a sphere or disk $\Delta$ \Defn{$B$-decomposable} if it is a simplex or there exists a vertex $w$ such that $\Delta-w$ is $B$-decomposable of the same dimension. We call it \Defn{$A$-decomposable} if the deletion is also $A$-decomposable, and ${\Lk}_w \partial \Delta$ is \Defn{$A$-decomposable} as well. We also point out the relation to work of Murai and Babson-Nevo: $A$-decomposable spheres are strongly edge-decomposable in the sense of \cite{BabsonNevo, MuraiHL}.

So, to summarize, we have understood how to prove the hard Lefschetz theorem for $A$-decomposable spheres. Though if we knew the Lefschetz theorem for ${\Lk}_w \partial \Delta$ for any other reason (say, induction on the dimension), $B$-decomposability would be enough.

\section{Pairings}

\subsection{Biased pairing}

We noted in the last example that the Lefschetz theorem is important for the induction. But we also noted that a property of the Poincar\'e pairing is important. And as we saw, these are in fact equivalent. Precisely, consider a triangulated $(d-1)$-ball or sphere $\Delta$ in $\R^d$. We say that it satisfies the \Defn{biased pairing property} (with respect to its boundary) if the map
\[\AR^k(\Delta, \partial \Delta)\ \longrightarrow\ \AR^k(\Delta)\]
is an injection for all degrees $k\le \frac{d}{2}$, or equivalently, that the pairing
\[\AR^k(\Delta,\partial \Delta)\ \times\ \AR^{d-k}(\Delta,\partial \Delta)\ \longrightarrow\ \AR^{d}(\Delta,\partial \Delta)\ \cong\ \R\]
is perfect. Note that this property is trivial for spheres.
If $k=\frac{d}{2}$, then this is equivalent to the middle Lefschetz isomorphism on the boundary of the disk (see \cite[Proposition 5.8]{AHL}), explaining the two different perspectives above. In other words, identifying what is necessary to do to show the non-degeneracy of the Poincar\'e product on kernel and image is easy in the case of $B$-decomposability.

\subsection{Biased Poincar\'e duality} This is a special case of {biased Poincar\'e duality} introduced in \cite[Section 5]{AHL}. This, in essence, has the goal of  extending the reach of the proof above: Instead of requiring that we look at disks in each step, we want to understand what happens in other cases to the pullbacks of kernel and image of the previous map. 

Let $\Sigma$ be a sphere of dimension $(d-1)$, with Poincar\'e pairing
\[\AR^k(\Sigma)\ \times\ \AR^{d-k}(\Sigma)\ \longrightarrow\ \R.\]
We say that an ideal $I$ in $\AR(\Sigma)$ satisfies \Defn{biased Poincar\'e duality}, if this pairing, restricted to $I$, is non-degenerate for $k\le \frac{d}{2}$. A thing to note:

If $\Delta$ is a disk in $\Sigma$ of the same dimension, and $\Gamma$ the complementary disk, then biased Poincar\'e duality for the ideal
\[\ker[\AR(\Sigma)\ \longrightarrow\ \AR(\Gamma)]\]
is the biased pairing property for $\Delta$. And as we shall see, that biased pairing property boils down to a Lefschetz theorem on $\partial \Gamma$.
Such ideals, obtained as kernels of
\[\AR(\Sigma)\ \longrightarrow\ \AR(X)\]
where $X$ is any subcomplex, are tremendously important. We also say therefore that biased Poincar\'e duality is satisfied at $X$ if and only if it is satisfied with respect to the associated ideal.

\subsection{Lefschetz and Hall-Laman} Let us finally come to the strongest notion: The biased pairing property is a weaker form of the hard Lefschetz property stated before, which we recall applies if with respect to some $\ell$ in $\R^1(\Delta)$ 
\[\AR^k(\Delta, \partial \Delta)\ \xrightarrow{\ \cdot \ell^{d-2k} \ }\ \AR^{d-k}(\Delta)\]
is an isomorphism for all $k\le \frac{d}{2}$. This map is naturally factored as
\[\AR^k(\Delta, \partial \Delta)\ \xrightarrow{\ \cdot \ell^{d-k} \ }\ \AR^{d-k}(\Delta,\partial \Delta)\ \longrightarrow\ \AR^{d-k}(\Delta).\]
This is a special case of the Hall-Laman relations of \cite{AHL}.

\section{Cones and boundaries}\label{sec:beymid}

Next, the proof of the hard Lefschetz theorem is reduced to the \Defn{middle case} when 
\[k\ =\ \frac{d}{2}.\] Note that in this case, both properties coincide.

To this end, assume we wish to prove the hard Lefschetz property for $\Delta$, which is a simplicial $(d-1)$-sphere or ball. Let $\cone \Delta$ denote the cone over $\Delta$ with apex $\mbf{n}$, realized in $\R^{d+1}$. Let $\uppi$ denote the projection along $\mbf{n}$, and let $\tet=h\cdot\x$, where $h$ is the vector of coordinates of the individual vertices in coordinate direction $\mbf{n}$. The cone lemmas then give:
	
	\begin{lem}[{\cite[Lemma 7.6]{AHL}}]\label{lem:midred}
Considering $\cone\Sigma$ realized in $\R^{d+1}$, and $k< \frac{d}{2}$, the following two are equivalent:
		\begin{compactenum}[$(1)$]
					\item The Hall-Laman relations for $\cone\Delta$
					with respect to $x_{\mbf{n}}$ and in degree $k+1$	.								
					\item The Hall-Laman relations for 
$ \uppi\Delta$			
with respect to $\tet$ in degree $k$.
\end{compactenum}
\end{lem}

This iteratively reduces to the middle case.

\begin{proof}
We may assume without loss of generality that $\tet=x_{\mbf{n}}$ in  $\AR^{\ast}(\cone\varDelta)$.
Consider then the diagram
\[\begin{tikzcd}[column sep=5em]
 \AR^{k}(\uppi\varDelta, \partial\uppi \varDelta) \arrow{r}{\ \ \ \ \cdot \tet^{d-2k}\ \ \ \  } \arrow{d}{\sim } & \AR^{d-k}(\uppi\varDelta) \arrow{d}{\sim } \\
\AR^{k+1}(\cone\varDelta, \partial \cone\varDelta) \arrow{r}{\ \ \ \cdot x_{\mbf{n}}^{d-2k-1}\ \ \ } & \AR^{d-k}(\cone\varDelta)
\end{tikzcd}
\]
Where the first vertical map is defined by the composition of cone lemmas
\[\AR^{k}(\uppi\varDelta, \uppi\partial \varDelta)\ \cong\ \AR^{k}(\cone\varDelta, \cone \partial \varDelta)\ \xrightarrow{\ \cdot x_{\mbf{n}} \ }\ \AR^{k+1}(\cone\varDelta, \partial \cone\varDelta).\]
and the second vertical map is simply the cone lemma. An isomorphism on the top is then equivalent to an isomorphism of the bottom map.
\end{proof}

Once we arrive there, another miracle happens, because now Lefschetz has become biased pairing:

\begin{lem}[{\cite[Lemma 5.6]{AHL}}]\label{lem:PLinv}
	A PL homeomorphism $\upvarphi:\Delta\rightarrow \Delta'$ of balls $\Delta,\ \Delta'$ in $\R^d$ that restricts to the identity on the boundary preserves the pairing property. That is, $\Delta$ satisfies the biased pairing property if and only if $\Delta'$ does.	
\end{lem}

This is a rather simple application of the decomposition theorem \cite{BBDG} in its simplest form: We want to show that passing from $\Delta$ to $\Delta'$ does not affect the non-degeneracy of the pairing
\[\AR^k(\Delta,\partial \Delta)\ \times\ \AR^{d-k}(\Delta,\partial \Delta)\ \longrightarrow\ \R\]
For this, we may restrict to single stellar subdivisions and inverses when passing from $\Delta$ to $\Delta'$, and the decomposition theorem shows that the contribution of such a modification splits off orthogonally.

In particular, this lemma allows us to arbitrarily refine, and therefore make simpler, the triangulation of $\Delta$ in the interior, only having to keep the boundary intact.

\section{So, how do we proceed?}\label{sec:beymid}

We now have reduced ourselves to prove the biased pairing property for a disk $\Delta$ of odd dimension $(2d-1)$, and specifically do so in dimension $d$. Equivalently, we can take the double of $\Delta$ (two copies of $\Delta$, identified along the boundary). In that sphere $\Sigma$, we are now expected to prove biased Poincar\'e duality in degree d and with respect to the ideal generated by $\partial \Delta$. Which in turn necessitates proving the middle Lefschetz isomorphism on the same. We are caught in a loop.

Or are we?

First, notice we only care about the $(d-1)$-skeleton, say $S$ of $\partial \Delta$. Biased Poincar\'e duality in degree $d$, with respect to the former or the latter, are equivalent, as the associated ideals coincide in degree $d$, because $\AR^d(\partial \Delta) \cong \AR^d(S)$.

Now, the following main idea comes into play: Say we find another hypersurface $E$, possibly with boundary, that is an envelope for $S$ in degree $d$, that is, $\AR^d(E) \cong \AR^d(S)$. If biased Poincar\'e duality is implied a Lefschetz statement on $E$. 

This is indeed the case:

Consider for this purpose the complement $\widetilde{\Sigma}$ of ${E}$ in $\Sigma$, where  ${E}$ is a $(d-2)$-acyclic rational hypersurface (with induced boundary) in $\Sigma$, and the double $\mathrm{D}{E}$ of ${E}$.

		The open manifold $\widetilde{\Sigma}$ can be compactified canonically to a compact manifold $\widetilde{\mr{D}}\Sigma$ with boundary $\mathrm{D}{E}$.
	
Let $\uppi$ denote the general position projection to a hyperplane $H$, and ${\tet}$ the height over that projection, so that
	\[\AR^k(X)\ = \ \bigslant{\AR^{k}(\uppi\hspace{0.3mm} X)}{ {\tet}\ \AR^{k-1}(\uppi\hspace{0.3mm} X)}\]
	for a complex $X$ in $\R^d$.
	
		\begin{lem}\label{lem:reductodouble}
In the situation above, assume that 
\begin{compactenum}[(1)]
		\item The map
		\[\AR^{d-1}(\uppi\hspace{0.3mm} {E},\partial \uppi\hspace{0.3mm} {E}) \xrightarrow{\ \cdot {\tet}\ }\ \AR^d(\uppi\hspace{0.3mm} {E}).\]
		is injective.
\item 		The map 
\[(H^{d-1})^{\binom{2d}{d}}({\widetilde{\mr{D}}\Mu})\ \longrightarrow\  \AR^d(\mathrm{D}{E}),\]
which factors as
\[\begin{tikzcd}
		(H^{d-1})^{\binom{2d}{d}}({\widetilde{\mr{D}}\Mu})\ \arrow{r}{}\arrow{dr}{}&  \AR^d(\mathrm{D}{E})\\
		&(H^{d-1})^{\binom{2d}{d}}({{\mr{D}}E}) \arrow{u}{}
		\end{tikzcd}
		\]
is an injection.
\end{compactenum}
Then $\Sigma$ satisfies biased Poincar\'e duality with respect to $E$.
	\end{lem}
	
The main work of \cite{AHL} is to show that both conditions are, in fact, Lefschetz type properties, and can be shown iteratively using the perturbation lemma. At least, provided $E$ is nicely decomposable.	Unfortunately, neither $A$ nor $B$-decomposability are enough.

\section{Conclusion}

One can show that every triangulated manifold $\Mu$, after some subdivision not involving its boundary, becomes \Defn{$C$-decomposable}, that is, there exists an order on a subset of the interior vertices of $\Mu'$ such that, for any initial segment $W$ and the next vertex $w$
\begin{compactenum}[(1)]
\item $\bigcup_{v\in W} \st_v \Mu'$ is a submanifold of $\Mu'$, and
\item for the final segment $\overline{W}$, $\bigcup_{v\in \overline{W}} \st_v \Mu'=\Mu'$, and
\item $\lk_w \Mu'\cap(\partial \Mu'\cup \bigcup_{v\in W} \st_v \Mu')$ and $\lk_w \bigcup_{v\in W} \st_v \Mu'$
 is $C$-decomposable of codimension one, or empty.
\end{compactenum}	
Finally, we call a set of vertices, seen as a $0$-dimensional manifold, $C$-decomposable by default. This notion is slightly weaker than what \cite{AHL} calls a m\'etro, where additional homological assumptions are made to simplify the proof.

We now have the following basic construction:

\begin{lem}Assume $S$ is a $(d-1)$-dimensional complex in a PL $(2d-1)$-sphere, or it is the $(d-1)$-skeleton of a $(2d-2)$-sphere $\Sigma'$, living in the suspension $\Sigma$ of the latter. Then, assuming the middle Lefschetz theorem for $(2d-4)$-spheres, there exists a hypersurface $E$ in a subdivision of $\Sigma$ that does not affect $S$ such that
\begin{compactenum}[(1)]
\item $E$ contains $S$, and is an envelope for $S$ in degree $d$.
\item $E$ is $C$-decomposable.
\end{compactenum}
\end{lem}

\begin{figure}[h!tb]
			\begin{center}
				\includegraphics[scale = 0.7]{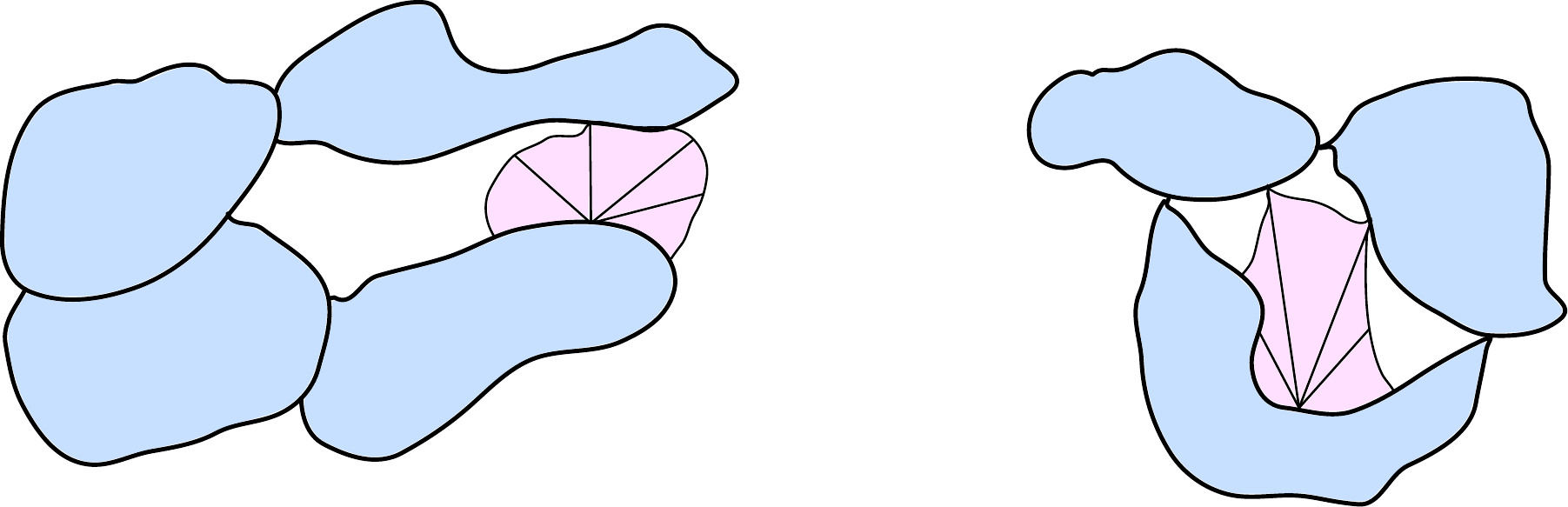}
				\caption{$C$-decomposability allows for simpler analysis of kernels and images in the pullback of a new vertex.}
				\label{decomp}
			\end{center}
\end{figure}

The key point of \cite{AHL} is that, while not as simple as the ideals that arise in the case of $B$-decomposability, $C$-decomposability guarantees they are simple enough to control to allow for the understanding in each application step of Lemma~\ref{lem:perturbation}. 

The observation is that the two conditions guarantee that in each step, the kernel of the map $``{\sum_{v\in W}}"x_v$, pulled back along $x_w$, can be analyzed once again using Lefschetz theorems on
$\lk_w \Mu'\cap(\partial \Mu'\cup \bigcup_{v\in W} \st_v \Mu')$. 

Indeed, this is even simpler to see and understand by working the orthogonal complement, that is, we need to understand the image of
$``{\sum_{v\in W}}"x_v$, intersected with the ideal $\langle x_w \rangle\in \AR^\ast(\Mu',\partial \Mu')$, thereby transferring us to $\AR^\ast(\st_w \Mu')$. 

But now observe the very best thing: This intersection vanishes under multiplication with $x_w$ as all the vertices of $W$ are at distance at least two from $w$ by the first condition. That condition then guarantees that the image is exactly consisting of the primitive elements under the multiplication with $x_w$ on $\lk_w \Mu'\cap(\partial \Mu'\cup \bigcup_{v\in W} \st_v \Mu')$. The hard Lefschetz property on the latter then implies that the ideal satisfies biased Poincar\'e duality in $\lk_w \Mu'$.
As this is, in turn, a ($C$-decomposable) manifold, it in turn closes the induction. 

\section{Via edge contractions}

The alternative proof uses, instead of the perturbation lemma, edge contractions. Edge contractions are a well-known and used to prove Lefschetz type theorems, but a major problem arises when trying to do this in full generality. As observed first by Whiteley \cite{Whiteley} (see also \cite{BabsonNevo, MuraiHL}, it is easy to ensure the Lefschetz property for a sphere $\varSigma$ if, after a contraction of the edge $e$ of $\varSigma$, the resulting complex $\varSigma'$ is (1) still a sphere, satisfies (2) the Lefschetz theorem and so does (3) the link of the contracted edge. Conditions (2) and (3) are easily satisfiable by induction, but it is condition (1) that breaks our neck. Indeed, it is not clear that we can find a contractible edge. And indeed, starting from dimension $3$, there exist spheres that are not the simplex, but still have no contractible edges (see again \cite{MuraiHL}).

The theory we introduced offers us a way out: We apply edge contraction to a suitable envelope of the $(d-1)$-skeleton instead.

Instead of showing, for an even dimensional sphere $\varSigma'$ (say, of dimension $2d-2$) the Lefschetz isomorphism between $\AR^{d-1}(\varSigma')$ and $\AR^{d}(\varSigma')$, we show the equivalent claim that in a suspension $\varSigma$ of the same, the $(d-1)$-skeleton induces an ideal that satisfies biased Poincar\'e duality.

But once again, we can try to find a sufficiently nice hypersurface envelope $E$ for that skeleton that is contractible. In other words, we call a hypersurface envelope $E$ for a $(d-1)$-complex $X$ then $D$-contractible if there exists an edge $e$ such that
\begin{compactenum}[(1)]
\item the contracted complex $E'$ is an envelope for the contracted complex $X'$ in degree $d$ and
\item $E'$ contains as a subcomplex a hypersurface $\widetilde{E}$ that contains $X'$ and is $D$-contractible, or $X'$ is of a dimension lower than $d-1$.
\end{compactenum}

With this, we can use Whiteley's induction and prove the theorem via edge contractions. We only need the following lemma:

\begin{lem}Assume $S$ is a $(d-1)$-dimensional complex in a PL $(2d-1)$-sphere, or it is the $(d-1)$-skeleton of a $(2d-2)$-sphere $\Sigma'$, living in the suspension $\Sigma$ of the latter. Then, assuming the middle Lefschetz theorem for $(2d-4)$-spheres, there exists a $D$-contractible hypersurface envelope $E$ for $S$ in a subdivision of $\Sigma$.
\end{lem}

\section{Open questions}

There are several open questions concerning the construction of Lefschetz elements, hard Lefschetz type theorems. Let us stay with the construction in detail: To prove the hard Lefschetz theorem for spheres, we needed to discuss more general manifolds:

\begin{qst}
Is it possible to prove the Lefschetz theorem for spheres without the detour over manifolds and hypersurfaces?
\end{qst}

Secondly, it remains to see whether there is a symmetric, or perhaps even homological version of the perturbation lemma, that does not rely on analyzing two vertices at a time.

\begin{qst}
Is it possible to combine more than two maps at once in Lemma~\ref{lem:perturbation}?
\end{qst}

Finally, we could ask for a characterization of "good" artinian reductions of face rings of spheres, such that a Lefschetz element can be found. From the above, we can distill the following way to find the right Artinian reductions:

\begin{thm}
Consider a $(2d)$-sphere $\Sigma$ realized in $\mathbb{R}^{2d+1}$. Assume that for every vertex $v$ of $\Sigma$, and every $(d-1)$-dimensional subcomplex of its link $S$, every embedding of $S$ into a $(2d-1)$-sphere satisfies biased Poincar\'e duality. Then $\Sigma$ has a Lefschetz element for the middle isomorphism.
\end{thm}

We draw a useful corollary.

\begin{cor}
If $\Sigma$ is a $2$-sphere in $\R^3$, then it has a Lefschetz element if no three vertices lie in a hyperplane through the origin.
\end{cor}

This leaves open the search for explicit coordinates for spheres to satisfy Lefschetz thoerems. To avoid cheating via transcendental extensions, we ask:
\begin{prb}
Is there a natural set of rational coordinates/rational Artinian reduction for a sphere $\Sigma$ to satisfy the Lefschetz theorem?
\end{prb}

	\textbf{Acknowledgements:} Karim Adiprasito was supported by the European Research Council under the European
Union's Seventh Framework Programme ERC Grant agreement ERC StG 716424 - CASe, the Israel Science Foundation under ISF Grant 1050/16 and the Knut and Alice Wallenberg foundation. The author is grateful to Giovanni Cerulli Irelli and Corrado De Concini, as well as the Sapienza University of Rome, for their hospitality, and the audience of the minicourse for many helpful comments, and to Johanna Steinmeyer for pointing out several typos, as well as Kalle Karu for useful comments.

%	{\small
%		\bibliographystyle{myamsalpha}
%		\bibliography{ref}}
%

\providecommand{\bysame}{\leavevmode\hbox to3em{\hrulefill}\thinspace}
\providecommand{\MR}{\relax\ifhmode\unskip\space\fi MR }
\providecommand{\MRhref}[2]{%
  \href{http://www.ams.org/mathscinet-getitem?mr=#1}{#2}
}
\providecommand{\href}[2]{#2}

\end{document}